\documentclass[12pt]{amsart}
\usepackage{amsmath,amscd,amsthm,amsfonts, amssymb,amsxtra}
\usepackage{epsf}
\usepackage[all]{xy}
\newtheorem{thm}{Theorem}[section]  
\newtheorem*{un-no-thm}{Theorem}
\newtheorem{lem}[thm]{Lemma}         

\newtheorem*{assertion}{Assertion}
\newtheorem{bigthm}{Theorem}

\theoremstyle{definition}

\theoremstyle{definition}
\theoremstyle{definition}

\theoremstyle{remark}
\newtheorem{rem}[thm]{Remark}

\newtheorem*{acks}{Acknowledgements}

\begin{document}
\title{On C.T.C.\ Wall's Suspension Theorem}
\date{\today}
\author{Mokhtar Aouina}
\address{Wayne State University, Detroit, MI 48202}
\email{aouina@math.wayne.edu}
\author{John R. Klein}
\address{Wayne State University, Detroit, MI 48202}
\email{klein@math.wayne.edu}
\begin{abstract} The concept of {\it thickening}
was systematically studied by C.T.C.\ Wall
in \cite{Wall4}. The {\it suspension theorem} of that paper
is an exact sequence relating the $n$-dimensional thickenings
of a finite complex to its $(n{+}1)$-dimensional ones. 
The object of this note is to fill in what we believe is
a missing argument in the proof of that theorem.
\end{abstract}
\thanks{The second author is partially supported by  NSF Grant DMS-0201695.}
\maketitle
\setlength{\parindent}{15pt}
\setlength{\parskip}{1pt plus 0pt minus 1pt}
\def\Top{\bold T\bold o \bold p}
\def\Sp{\bold S\bold p}
\def\vo{\varOmega}
\def\vs{\varSigma}
\def\smsh{\wedge}
\def\flush{\flushpar}
\def\id{\text{id}}
\def\dbslash{/\!\! /}
\def\codim{\text{\rm codim\,}}
\def\:{\colon}
\def\holim{\text{\rm holim\,}}
\def\hocolim{\text{\rm hocolim\,}}
\def\hodim{\text{\rm hodim\,}}
\def\hocodim{\text{hocodim\,}}
\def\Bbb{\mathbb}
\def\bold{\mathbf}
\def\Aut{\text{\rm Aut}}
\def\cal{\mathcal}
\def\frak{\mathfrak}

\section{Introduction} 

One of the classical industries in differential topology is to 
list up to diffeomorphism 
all manifolds within a given homotopy type. 
In dimensions $\ge 6$, assuming that the
homotopy type of the boundary is held fixed, the
surgery exact sequence \cite{Wall_2nd} reduces the problem to
bundle theory and the stable algebra of bilinear forms over the integral
group ring of the fundamental group. 

However, if the homotopy type of the
boundary is not held fixed, 
the methods of surgery theory do not apply. 
Nearly forty years ago, C.T.C.\  Wall  \cite{Wall4} 
developed tools 
to study this version of the problem.
Of these, probably the most computationally 
useful is the {\it suspension theorem}
(cf.\ \cite[\S5]{Wall4}), which is a variant of
the classical EHP sequence of I.M.\ James \cite{James_triad}. 
This note is a comment on the suspension theorem.

\subsection*{Thickenings} 

Let $K$ be a connected finite CW complex of dimension $\le k$.
An {\it $n$-thickening} of $K$ consists of a simple homotopy equivalence
$$
f\: K \to M
$$
in which $M$ is a compact smooth manifold of dimension $n$ such that
the inclusion $\partial M \to M$ is $1$-connected and induces an
isomorphism of fundamental groups (such manifolds are said to satisfy
the {\it $\pi$-$\pi$ condition}). We will assume that
$k \le n{-}3$ and $n \ge 6$.

Two $n$-thickenings $f_0\: K \to M_0$ and $f_1\:K \to M_1$ are said
to be {\it equivalent} if there is a diffeomorphism $h\:M_0 \to M_1$
such that $hf_0$ is homotopic to $f_1$. This gives an equivalence
relation on the set of $n$-thickenings; let 
$$
{\cal T}_n(K)
$$
denote the associated set of equivalence classes.

\subsection*{Suspension}
The operation obtained by taking the cartesian product with the unit
interval gives a map
$$
E\:{\cal T}_n(K) \to {\cal T}_{n{+}1}(K)
$$
called {\it suspension.} Let us say that an $(n{+}1)$-thickening
{\it compresses} if its associated equivalence class is in the
image of $E$.

Wall's {\it suspension theorem} provides criteria for
deciding when a thickening compresses.
Wall assumes
\begin{itemize}
\item $K$ is $(2k{-}n)$-connected, and
\item $3k {+}2 \le 2n$.
\end{itemize}
He then constructs a function
$$
H\:{\cal T}_{n{+}1}(K) \to [K {\times} K,S^{n+1}]
$$
where the target is the homotopy classes of maps from $K {\times} K$
to $S^{n{+}1}$ (we will review the definition of $H$ below).
Wall then shows that a thickening compresses if and only if
it vanishes under $H$. This is the content of one half of the 
suspension theorem. The other half gives criteria for deciding the
extent to which $E$ is injective (we will not be concerned with the
latter in this paper).

\subsection*{Definition of $H$}
Let $f\: K \to M$ be an $(n{+}1)$-thickening representing an element
$x \in {\cal T}_{n{+}1}(K)$. Let $C$ be the complement
of a tubular neighborhood of the diagonal embedding $M \to M \times M$. Then
we have an identification
$$(M\times M)/C \,\,  \cong \,\, M^\tau\, ,
$$ 
where 
$M^\tau$ denotes the Thom space tangent
bundle $\tau$ of $M$. 
Choose a basepoint for $K$. Applying $f$, we obtain a basepoint for $M$. Identify $S^{n{+}1}$ with the one point compactification of
the fiber of $\tau$ at the basepoint. This gives an inclusion 
$$
S^{n{+}1} \subset M^\tau\, .
$$
Wall observes that the range assumptions imply that this
inclusion is $(2k{+}1)$-connected.

The composite map 
$$
\begin{CD}
K \times K @> f\times f >> M \times M @>\text{collapse } C>> 
(M \times M)/C  = M^\tau
\end{CD}
$$
Gives rise to a homotopy class in $[K \times K,M^\tau]$.
Since $K \times K$ has dimension $\le 2k$, the induced map
of homotopy sets
$$
[K\times K,S^{n{+}1}] \to [K\times K,M^\tau]
$$
is a bijection. Wall defines 
$$H(x) \in [K\times K,S^{n{+}1}]$$ to be
the unique element that pushes forward to the above
homotopy class. 

\subsection*{The key step in Wall's proof}
Wall wishes to show that the triviality of
$$
H(x) \in [K\times K,S^{n{+}1}]
$$
implies that the inclusion $\partial M \to M$ admits a section up to homotopy.

If there is such a homotopy section $s\:M \to \partial M$, then 
by the numerical assumptions and the Stallings-Wall embedding
theorem (\cite{Stallings_embed},\cite{Wall4}), the map
$$
\begin{CD}
K @> f >> M @> s >> \partial M
\end{CD}
$$
admits an {\it embedding up to homotopy,} i.e., 
there is a compact codimension zero submanifold $N \subset \partial M$
and a simple homotopy equivalence $g\:K \to N$ such that the composite
$$
\begin{CD}
K @> g >> N @> \subset >> M
\end{CD}
$$
is homotopic to $f$. Then, using the
$s$-cobordism theorem, one sees that
$g\:K \to N$ is a representative of a compression of $x$.

In summary, the proof hinges on the following:

\begin{assertion}{\rm (Wall).}  In the given range, $H(x)$ vanishes
if and only if $\partial M \to M$ admits a homotopy section.
\end{assertion}

On page 87 of \cite{Wall4}, a proof of this assertion is purportedly
given. We cite the relevant paragraph here:

\begin{quote} \quad Now let $M'$ be the complement of a collar neighborhood 
of $\partial M$ in $M$, $i \: M' \to M$ the inclusion. Since $K,M$, and $M'$
all have the same homotopy type, the assertion above shows that $f\times i\:
K \times M' \to M \times M$ is homotopic to a map into $C$, and hence to a
map into $\text{\rm Int}(M\times M - \Delta M)$. But the projection
onto the second factor $\text{\rm Int}(M)$ is a fibration (it is well known
to be locally trivial); the above homotopy projects to one in 
$\text{\rm Int}(M)$, whose inverse lifts to a homotopy of the constructed
map $K \times M' \to C$ to map of the form $g' \times i$. Thus
$f \simeq g'$, and the image of $g'$ lies in a collar neighborhood of
$\partial M$; a further evident homotopy  sends $g'$ to a map $g\: K \to
\partial M$.
\end{quote}

{\flushleft (note: we have changed Wall's 
notation to conform with ours, but with this exception, what is displayed here 
is verbatim).}  In line 3 of the citation, 
the expression `the assertion above' amounts to
the statement that the cofibration sequence
$$
C \to M \times M \to M^\tau
$$
is a homotopy fiber sequence in dimensions $\le 2k$. This
follows from the Blakers-Massey theorem and the range assumption.
Hence the first part of the argument is clear to us.

However, the latter half of the paragraph eludes us. 
Specifically, we do not understand how the map $g' \times i$ is produced having
the required property: namely, the map $g'\: K \to M$ 
is supposed to satisfy $g'(x) \ne i(y)$ for all 
$x \in K$ and $y \in M'$. It seems to us
that the map $g'$ obtained from Wall's argument really depends on
both factors, and so is a map of the form  $g' \:K \times M' \to M$.
Unfortunately a map of the latter kind is 
insufficient for producing the desired map $g$.

Furthermore, once $f \times i$ is homotoped to a map into
$C$, {\it the range assumption is not used anymore in the argument.}
Since $f$ and $i$ are homotopy equivalences, one can
reformulate what Wall writes as saying that, without
any range restrictions, the only obstruction to
finding a homotopy section of the inclusion $\partial M \to M$
is to find a homotopy factorization of the identity map of
$M \times M$ through $M \times M - \Delta M$. We think this
is unlikely in general to be true, but we do not have a counterexample.

Rather than trying to sort out Wall's argument, we will provide
one of our own.

\subsection*{Reformulation}
The above description of $H$ used
the collapse map
$$
M \times M \to M^\tau
$$
which make into a based map 
$$
c_0\: M_+ \smsh M_+ \to M^{\tau} \, ,
$$
where $M_+$ is the effect of adjoining a disjoint basepoint to $M$
and we have used the identification $(M{\times} M)_+ = M_+ \smsh M_+$.

It is clear from the range assumption and the definition of $H$
that $H(x)$ vanishes if and only if $c$ is stably null homotopic.

We are ready to formulate the main result of this paper.

\begin{bigthm} \label{wall_fix} 
Let $M$ be a compact manifold of dimension $n{+}1$
satisfying the $\pi$-$\pi$ condition. Assume that $M$ has the
homotopy type of a CW complex of dimension $\le k$. Furthermore,
assume
\begin{itemize}
\item $M$ is $(2k{-}n)$-connected, and
\item $3k {+}1 \le 2n$.
\end{itemize}
Then the inclusion $\partial M \to M$ admits a section
up to homotopy if and only $c_0\: M_+ \smsh M_+ \to M^\tau$
is stably null homotopic.
\end{bigthm}

Note that the range of our theorem is one dimension better than
the one used by Wall.





\begin{acks} This paper arose out of discussions related to
the first author's Ph.D.\ thesis. We are indebted to Bill
Richter for suggesting that a proof of the main
result ought to exist using some form of duality.
We are also greatful to Bruce Williams and Michael Weiss for help
in tracking down references to theorem \ref{pd}. 
\end{acks}

\section{The proof of Theorem \ref{wall_fix}}

We first observe that the inclusion  
$$
M \times \partial M \subset M \times M
$$
is homotopic to a map into $C$ (to see this, let $M_0$ be the complement of
a collar neighborhood of $\partial M$. Then $\partial M \times M_0$ is
a subspace of $C$ and $M_0 \subset M$ is a homotopy equivalence). This means
that the map $c_0$ factors up to homotopy as
$$
\begin{CD}
M_+ \smsh M_+ @>\text{id} \smsh p >>  M_+ \smsh M/\partial M @>c >> M^\tau
\end{CD}
$$
in which $p \: M_+ \to M/\partial M$ is the evident quotient map and the map
$c$ is the map on quotients induced by the map of pairs
$$
(M \times M,M_0  \times \partial M) \to (M\times M,C) 
$$
(here we are implicitly identifying 
the quotient $(M\times M)/(M_0 \times \partial M)$ with
the quotient $M_+ \smsh M/\partial M = (M\times M)/(M \times \partial M)$).

The map $c$ has an adjoint (cf.\ the exponential law)
$$
c^\# \: M/\partial M \to \text{maps}(M_+,M^\tau)\,.
$$

\begin{lem} \label{c_sharp} The map $c^\#$ is $(k +1)$-connected.
\end{lem}

 We will sketch a proof of \ref{c_sharp} in \S3.

\begin{proof}[Proof of Theorem \ref{wall_fix}, assuming \ref{c_sharp}]
As $M$ has the homotopy type of a
CW complex of dimension $\le k$, it follows from \ref{c_sharp} that
the map $p\: M_+ \to M/\partial M$ is null homotopic if and only if 
the map $$c^\sharp \circ p\:M_+ \to \text{maps}(M_+,M^\tau)$$ 
is null homotopic.

The $\pi$-$\pi$ condition,
the relative Hurewicz theorem and Poincar\'e duality, imply that
the map $\partial M \to M$ is $(n{-}k)$-connected.
From this it follows that $\partial M$ is $\min(n-k-1,2k-n)$-connected
since $M$ is $(2k-n)$-connected. As $2k - n \le n - k - 1 $ if and only
if $3k + 1 \le 2n$ (which is part of our assumptions), the minimum coincides
with $2k-n$. Consequently, $\partial M$ is also $(2k-n)$-connected.

By the Blakers-Massey theorem, the cofiber sequence
$$
\partial M \to M \to M/\partial M
$$
is a fibration up to homotopy up through dimension $k$. This means
that the map
$$
\partial M \to \text{hofiber}(M \to M/\partial M)
$$
is $k$-connected. As $M$ has the homotopy type of a CW complex of dimension
$\le k$, we infer from this that $\partial M \to M$ has a homotopy section
if and only if the map $p\:M_+ \to M/\partial M$ is null homotopic.

By Lemma \ref{c_sharp}, $p$ is null homotopic if and only if $c^\# \circ p$ is.
But the adjoint of the latter is a map $M_+ \smsh M_+ \to M^\tau$,
which is easily checked to coincide with Wall's map $c_0$.
\end{proof} 

\section{The proof of Lemma \ref{c_sharp}}
  
The proof will involve a generalized version of Poincar\'e duality.

Let
$$
E(\tau) \to M
$$
be the spherical fibration whose fiber $E(\tau)_x$ over $x\in M$ is the 
one point compactified tangent space at $x$
(after a Riemannian metric is chosen,
we have an identification $E(\tau)_x \cong D(\tau)_x/S(\tau)_x$,
where $D(\tau)_x$ is the unit tangent disk and $S(\tau)_x$
is the unit tangent sphere at $x$).
This sphere bundle
comes equipped with a section ``at infinity,'' which we regard as a
basepoint for its space of sections  $\text{sec}(E(\tau) \to M)$.

One defines a based map
$$
d\:M/\partial M \to \text{sec}(E(\tau) \to M)\, ,
$$
in the following way: the 
tubular neighborhood theorem gives a homotopy pushout square
$$
\SelectTips{cm}{}
\xymatrix{
S(\tau) \ar[r]\ar[d] & M \times M - \Delta_M \ar[d]\\
D(\tau) \ar[r] & M \times M \, .
}
$$

We take the fiber
at $x \in M$ 
of this diagram along the second factor projection $M \times M \to M$.
This gives a homotopy pushout
$$
\SelectTips{cm}{}
\xymatrix{
S(\tau)_x \ar[r]\ar[d]& M - x  \ar[d]\\
D(\tau)_x \ar[r]& M \times x\, .
} 
$$
We therefore
have a preferred identification
$$
E(\tau)_x \,\, \simeq \,\,    \text{cone}(M - x  \to M \times x)\, \, 
$$
where the right side denotes the mapping cone of 
$M -x \to M \times x$. So we 
get a (collapse) map
$$
M \times x \to \text{cone}(M - x  \to M \times x) \simeq  E(\tau)_x \, .
$$
Adjointly, what we have so far produced is a map
$$
M \to  \text{sec}(E(\tau) \to M)\, .
$$

Let $M_0$ denote the complement of a collar neighborhood of the
inclusion $\partial M \subset M$. By construction, the composite
$$
\begin{CD}
M @>>> \text{sec}(E(\tau) \to M) @>\text{restrict}> \simeq >
 \text{sec}(E(\tau)_{|M_0} \to M_0)
\end{CD}
$$
has the property that it maps $\partial M$ to the basepoint of
the section space (this is because $\partial M$ and $M_0$ 
are disjoint). So, up to the above identifications, what
we have really produced is a map
$$
d\:M/\partial M \,\, \to\,\,  \text{sec}(E(\tau) \to M)\, .
$$

Let us next think of the function space
$\text{maps}(M_+, M^{\tau})$ as a section space
of the trivial fibration $$M^{\tau}\times M \to M\, .$$ 
Change
the notation of the latter to $\bar E(\tau) \to M$. By construction,
the map $c^\#$ factors up to homotopy as
$$
\begin{CD}
M/\partial M @> d >>  
\text{sec}(E(\tau) \to  M) @>q>>  \text{sec}(\bar E(\tau) \to  M)
\end{CD}
$$
in which $q$ is the map coming from the map of fibrations
$E(\tau) \to \bar E(\tau)$ defined using the evident inclusion of each fiber
$E(\tau)_x \to M^{\tau} \times x$.

\begin{lem}\label{q} The map $q$ is $(k{+}1)$-connected.
\end{lem}

\begin{proof} As noted in the introduction, the map of
fibers
$$
E(\tau)_x \to  M^\tau \times x
$$
is $(2k{+}1)$-connected. Obstruction theory then
shows the connectivity of $q$ is this 
number minus the dimension of $M$.
\end{proof}

\begin{lem} \label{d} The map $d$ is $(2n-2k+1)$-connected. 
\end{lem}

Before proving \ref{d}, we show why these results imply \ref{c_sharp}.
Observe $k + 1 \le 2n - 2k + 1$ is equivalent to $3k \le 2n$, so \ref{d} 
and our range assumption will
imply that $d$ is $(k{+}1)$-connected. 
This, together with \ref{q} and
the factorization $c^\sharp = q\circ d$ implies that $c^\#$ is 
$(k{+}1)$-connected, yielding \ref{c_sharp}.
\medskip

\begin{proof}[Proof of \ref{d}]
There is a commutative diagram
$$
\SelectTips{cm}{}
\xymatrix{ 
M/\partial M \ar[r]^{d\qquad} \ar[d] &  \text{sec}(E(\tau) \to M)\ar[d]\\
Q(M/\partial M) \ar[r]^{d^{\rm st}\quad } &  \text{sec}^{\text{st}}(E(\tau)\to  M)
}
$$
where $Q$ is the stable homotopy functor, $\text{sec}^{\text{st}}$ means
{\it stable sections} (= sections of the fibration whose fiber at $x \in M$
is $Q(E(\tau)_x)$\, ), and the vertical maps are the natural inclusions.
The bottom map $d^{\rm st}$ is a stable
version of $d$; one way to
see this is to  think of $Q(M/\partial M)$ as
a limit over $j$ of $\Omega^j (M\times D^j)/\partial (M\times D^j)$,
so we can take a version of the top horizontal map using
$M \times D^j$ instead of $M$ and then loop it $j$ times.
The details will be left
to the reader.
\medskip 

Since $M/\partial M$ is $(n{-}k)$-connected, the Freudenthal theorem
shows that the left vertical map is $(2(n-k)+1)$-connected. Since
$E(\tau)_x$ is $n$-connected, the map $E(\tau)_x \to Q(E(\tau)_x)$
is $(2n{+}1)$-connected. By elementary obstruction theory, we see that the
right vertical map is $(2n{-}k{+}1)$-connected. 
Using these observations,
Lemma \ref{d} is then a direct consequence of the following:

\begin{thm}[Poincar\'e Duality] \label{pd} For any
compact manifold smooth manifold $M$,
The map $d^{\rm st}$ is a weak equivalence (of infinite loop spaces).
\end{thm}
\end{proof}

\begin{rem} 
The earliest reference we could find for \ref{pd}
is in a paper of Dax \cite[Prop.\ 5.1]{Dax}, who formulates
it in terms of normal bordism theory, and proves
it using transversality. Our map $d^{\rm st}$
corresponds, via the Thom-Pontryagin
map, to the isomorphism denoted by $U$ in his proof.

In a 
paper of Weiss and Williams \cite[Prop.\ 2.4]{WWI},
a statement of the above kind
is proved using a map going
in the other direction (their statement allows 
for any spectrum of 
coefficients, not just the sphere spectrum).
Williams
has informed us that our map is not a homotopy inverse to the Weiss-Williams
map; one must first twist by a certain involution.

Another  proof of \ref{pd} which is more generally valid for 
Poincar\'e duality spaces (using a an alternative description 
of $d^{\rm st}$ as a `norm map') 
can be found in \cite[Thms.\ A,D]{Klein6}.

\end{rem}

\end{document}